\documentclass[12pt]{article}
\pdfoutput=1
\usepackage{amssymb}
\usepackage{amsmath}
\usepackage{amsthm}
\usepackage{color}
\usepackage{graphicx}
\usepackage{enumerate}
\usepackage{subfigure}

\newtheorem{thm}                   {Theorem}  [section]

\newtheorem{lem}          [thm]{Lemma}
\newtheorem{cor}           [thm]{Corollary}
\newtheorem{conj}         [thm]{Conjecture}

\theoremstyle{definition} \newtheorem{definition}{Definition}
\theoremstyle{definition} \newtheorem{remark}     [thm]{Remark}   
\theoremstyle{definition}    

\newcommand{\DCD}{\mathit{DCD}}

\usepackage{hyperref}

\begin{document}

\title{Danzer's configuration revisited}

\author{
Marko Boben \\
Faculty of Computer Science\\
University of Ljubljana\\
1000 Ljubljana, Slovenia\\ and\\
University of Primorska, IAM\\
Muzejski trg 2\\
6000 Koper, Slovenia\\
{\tt Marko.Boben@fri.uni-lj.si}\\
\and 
G\'abor G\'evay \\
Bolyai Institute\\ 
University of Szeged\\
Aradi v\'ertan\'uk tere 1.\\
6720 Szeged, Hungary\\
{\tt gevay@math.u-szeged.hu}\\
\and
Toma\v z Pisanski\thanks{The first and third author were supported
in part by the ARRS of Slovenia, grants: P1-0294 and N1-0011: GReGAS, and by the European Science Foundation. 
The second and third author were supported by the European Union and co-funded by the European Social Fund; 
project title: ``Telemedicine-focused research activities on the field of Mathematics, Informatics and Medical Sciences";
project number: T\'AMOP-4.2.2.A-11/1/KONV-2012-0073.
} \\
Faculty of Mathematics and Physics\\
University of Ljubljana\\
1111 Ljubljana, Slovenia\\
and\\
University of Primorska, FAMNIT\\
Glagolja\v{s}ka 8\\
6000 Koper, Slovenia\\
{\tt Tomaz.Pisanski@fmf.uni-lj.si}\\
}

\date{\today}
\maketitle

\newpage

\centerline{\emph{Dedicated to the memory of Ludwig Danzer (1927-2011)}}
\bigskip

\begin{abstract}
\noindent
We revisit the configuration of Danzer $\DCD(4)$, a great inspiration for our work. This configuration of type 
$(35_4)$ falls into an infinite series of geometric point-line configurations $\DCD(n)$. Each $\DCD(n)$ is
characterized combinatorially by having the Kronecker cover over the Odd graph $O_n$ as its Levi graph. 
Danzer's configuration is deeply rooted in Pascal's Hexagrammum Mysticum. Although the combinatorial configuration is highly 
symmetric, we conjecture that there are no geometric point-line realizations with 7- or 5-fold rotational symmetry;
on the other hand, we found a point-circle realization having the symmetry group $D_7$, the dihedral group of order 14. 
\end{abstract}

\noindent
{\small 
\medskip

\noindent
\emph{Math.\ Subj.\ Class.:  51A20, 52C30, 52C35, 05B30.}}

\bigskip

\section{Introduction} \label{Intro}

There is a remarkable point-line configuration of type $(35_4)$, due to Ludwig Danzer, which he never published~\cite{Gru08, GR90}. 
It is the first known geometric $(n_4)$ configuration with $n$ not divisible by 3. It has many interesting combinatorial, geometric and 
graph-theoretic relationships, and our goal in this paper is to present and discuss some of them. To put this object into a proper 
historical context, we should go back to Cayley, 1846. Cayley has given the following construction for Desargues' $(10_3)$ 
configuration~\cite{Cox50}. Let $P_1,\dots,P_5$ be five points in projective 3-space in general position, i.e.\ such that no four of 
them lie in a common plane. Then the 10 lines $P_iP_j$ and the 10 planes $P_iP_jP_k$ form a line-plane configuration such that 
its intersection with a plain containing none of the points  $P_i$ yields the point-line $(10_3)$ configuration of Desargues. Take now 
the dual construction. Instead of points, in this case we start from five planes in general position (i.e.\ no more than two planes meet 
in a line). Then every three of them meet in a point, and every two of them meet in a line, forming directly the $(10_3)$ point-line 
Desargues configuration.

Danzer's construction is completely analogous. It uses seven 3-spaces in general position in the 4-dimensional projective space. They meet 
by fours and by threes in 35 points and 35 lines, respectively; hence we obtain a $(35_4)$ configuration. (By suitable projection it can be 
carried over to a corresponding planar configuration). 

The possibility of a common generalization of Cayley's and Danzer's construction is now clear. Consider an $n$-dimensional projective space 
($n\ge 2$). In this space, take $2n-1$ hyperplanes in general position; this means that no more than $n$ of them meet in a point. We shall call 
such an arragement of hyperplanes in a projective $n$-space a {\em Cayley-Danzer arrangement\/}. Now $n$ hyperplanes of this arrangement 
meet in a point, and $n-1$ hyperplanes meet in a line. Thus we have altogether $\binom{2n-1}{n}$ points and  $\binom{2n-1}{n-1}$ lines 
(hence their number is the same). These points and lines can be considered as labelled with the $n$- and $(n-1)$-element subsets, respectively, 
of  a $(2n-1)$-set. Moreover, the incidence between the points and lines is determined by containment between the corresponding subsets. 
Hence it is clear that each point is incident with $n$ lines, and vice versa. Thus, after projection to a suitable plane, we have a planar 
point-line configuration of type
\begin{equation} \label{DanzerType}
\left( {\binom{2n-1}{n-1}\!}_n\, \right).
\end{equation}
Thus we see that both Desargues' $(10_3)$ and Danzer's $(35_4)$ configuration are particular cases of (\ref{DanzerType}); therefore, 
the class of configurations described above may be termed a {\em Desargues-Cayley-Danzer configuration\/}. We shall denote it by 
$\DCD(n)$.

In this paper we first discuss further interesting historical connections of Danzer's configuration (Section~\ref{Pascal}). Namely, we 
show that in fact it was already known in the nineteenth century, as a member of the infinite family of configurations associated with 
Pascal's {\em Hexagrammum Mysticum\/}. Then, in Section~\ref{Odd}, we explore some graph-theoretic connections of the configurations
$\DCD(n)$. In particular, we show the each $\DCD(n)$ is characterized combinatorially by having the Kronecker cover over the Odd 
graph $O_n$ as its Levi graph. In Section~\ref{Decomp} we present a decomposition of $\DCD(n)$. In addition to the subconfigurations 
found in this way, in Section~\ref{Coxeter} we establish the existence of another subconfiguration which is closely related to the 
well-known Coxeter graph. All these subconfigurations can be used to find explicit constructions for a graphical representation of
Danzer's configuration. Section~\ref{Constr} is devoted to this task. Finally, in Section~\ref {SectionCircles} we show that all the 
configurations studied in the former sections have not only point-line, but also point-circle realizations. In particular, we show
that although $\DCD(4)$ seems to admit no point-line realization with 7-fold symmetry, there are point-circle examples realizing
$D_7$ symmetry.

We recall some standard definitions. A {\em combinatorial\/} (or {\em abstract\/}) {\em configuration\/} of {\em type\/} $(p_q, n_k)$ is an incidence structure 
with sets $\mathcal P$ and $\mathcal B$ of objects, called points and blocks, such that
\begin{enumerate}[(1)]

\item
$|\mathcal P | = p$;

\item
$|\mathcal B | = n$;

\item
each point is incident with $q$ blocks;

\item
each block is incident with $k$ points;

\item
two distinct points are incident with at most one block.

\end{enumerate}
A  {\em point-line configuration\/} is a geometric incidence structure consisting of points and (straight) lines, in the simplest case in Euclidean 
and real projective plane, such that the conditions (1)--(4) hold. Note that in this case $(5)$ is fulfilled automatically. Instead of plane, the ambient 
space can be some higher-dimensional (Euclidean, or projective) space as well~\cite{Cox50, GG14}. Also, the set of lines can be replaced by 
a set of circles; in this case we speak of a {\em point-circle configuration\/}. Usually, it is defined in Euclidean plane (however, for more detailed 
investigation, it is appropriate to use inversive plane~\cite{GP}). Note that for a point-circle configuration condition $(5)$ does not necessarily 
hold; in case it still holds, we say that the configuration is {\em lineal\/}~\cite{GP}. 

It may happen that a point-line configuration $\mathcal C_1$ and a point-circle  configuration $\mathcal C_2$ are {\em isomorphic\/}; in 
other words, the underlying combinatorial configuration is essentially the same. In this case we say that $\mathcal C_1$ and $\mathcal C_2$ 
are different {\em geometric realizations\/} of the same abstract configuration.

Given a configuration of type $(p_q, n_k)$, if $p=n$, then $q=k$, and in this case the more concise notation $(n_k)$ is used. Formerly, such a 
configuration was called {\em symmetric\/}; however, it is appropriate to reserve this term for a configuration which has non-trivial automorphism, 
hence, following Gr\"unbaum~\cite{Gru09}, we use the term {\em balanced configuration\/} instead.

For further definitions and background material concerning configurations in general, the reader is referred to the recent monographs~\cite{Gru09, PS}.

\section[The 35 4 configuration associated with the complete Pascal hexagon]{The $\boldsymbol{(35_4)}$ configuration associated with the com\-plete Pascal hexagon} \label{Pascal}

Pascal's famous theorem states that if a hexagon is inscribed in a conic, then the three pairs of opposite sides meet in three points on 
a straight line. This line is called the {\em Pascal line\/}. Taking the permutations of six points on a conic, one obtains 60 different 
hexagons. Thus, the so-called {\em complete Pascal hexagon\/} determines altogether 60 Pascal lines. Steiner started to investigate 
the 60 Pascal lines, and he found that these lines meet in triples in 20 points (1828~\cite{Klu, Sal, CR}). These latter, which are now 
called {\em Steiner points\/}, lie, in turn, in quadruples on 15 {\em Pl\"ucker lines\/} (Pl\"ucker, 1929, see~\cite{CR}, and 1830, 
see~\cite{Klu}). Thus we get a $(20_3, 15_4)$ configuration, which we call the {\em Steiner-Pl\"ucker\/ configuration\/}.  Here 
we have the following incidence theorem (see Exercise II.16.6 in~\cite{VY} and Exercise 2.3.2 in~\cite{Cox74}).

\begin{thm} \label{TripleDesThm}
If three triangles are perspective from the same point, the three axes of perspectivity of the three pairs of triangles are concurrent.
\end{thm}

By Hesse, the $(20_3, 15_4)$ Steiner-Pl\"ucker configuration coincides with the configuration formed by the points and lines 
in this theorem (1850, see~\cite{Klu}). (For an illustration, see Figure \ref{SP_Labels} in Section~\ref{Decomp}.)

On the other hand, Kirkman discovered that the 60 Pascal lines intersect in threes in 60 additional points; in fact, the Pascal lines 
together with these \emph{Kirkman points} form a $(60_3)$ configuration (1849~\cite{Klu, Lev,CR}). Cayley  showed that the 
Kirkman points lie in triples on 20 new lines (1849~\cite{Klu,CR}). These 20 \emph{Cayley} lines were found to be concurrent 
in triples in 15 \emph{Salmon points} (Salmon, 1849, see~\cite{Klu,CR}). The 20 Cayley lines and the 15 Salmon points form a 
$(15_4, 20_3)$ configuration. Here we have again an incidence theorem, which is dual to the former (\cite{VY}, Exercise II.16.6)

\begin{thm} \label{TripleDesDualThm}
If three triangles are perspective from the same line, the three centres of perspectivity of the three pairs of triangles are collinear.
\end{thm}

The \emph{Cayley-Salmon configuration} formed by the 20 Cayley lines and the 15 Salmon points is the same as that obtained from the 
points and lines of this theorem. Thus the Steiner-Pl\"ucker configuration and this latter configuration are dual to each other. It is to be 
emphasized, however, that in spite of this dual correspondence, there is no polarity (in contrast to Hesse, 1868, see~\cite{Klu}), and not 
even a general projective correlation (in contrast to Schr\"oter, 1876, see~\cite{Klu}), which would carry the Pascal lines, Steiner points 
and Pl\"ucker lines into the Kirkman points, Cayley lines and Salmon points, respectively (cf.\ Remark~\ref{NoExtension} later). This 
was already known to Klug, who published a comprehensive work in 1898 on the configurations associated with the complete 
Pascal hexagon~\cite{Klu}.

Now, on page 34 of his book, Klug makes the following observation.

\begin{thm} \label{KlugObservation}
The  $(20_3, 15_4)$ Steiner-Pl\"ucker configuration and the $(15_4, 20_3)$ Cayley-Salmon configuration together form a $(35_4)$ configuration.
The 35 points and the 35 lines of this configuration are the 20 Steiner points and 15 Salmon points, and the 15 Pl\"ucker lines and 20 Cayley lines,
respectively. On each Pl\"ucker line there are four Steiner points, and on each Cayley line there are three Salmon points and one Steiner point; 
moreover, three Pl\"ucker lines and one Cayley line pass through each Steiner point, and four Cayley lines pass through each Salmon point.
\end{thm}

He also remarks that this configuration is nothing else than the planar projection of a figure formed
\begin{itemize}
\item by the vertices and edges of three tetrahedra, which are perspective from the same point; 
\item by the points and lines of intersection of (the projective hulls of) the corresponding edges and faces, respectively; 
\item and finally, by the projecting lines and by the lines of intersections of planes of perspectivity.
\end{itemize}

In this latter observation, the \emph{plane of perspectivity} is an analogue of one dimension higher of the axis of perspectivity, as 
it occurs in the following theorem on perspective tetrahedra (\cite{VY}, \S 17, Theorem 2; this theorem is already mentioned in one 
of Felix Klein's works, 1870, see Carver~\cite{Car}).

\begin{thm} \label{PerspectiveTetrahedra}
If two tetrahedra are perspective from a point, the six pairs of lines of the corresponding edges intersect in coplanar points, 
and the planes of the four pairs of faces intersect in coplanar lines; i.e.\ the tetrahedra are perspective from a plane.
\end{thm}

We remark that the configurations associated to the complete Pascal hexagon is a vast topic. In fact, there are infinitely many 
such configurations~\cite{Klu}; we do not pursue this topic further, since here we are only interested in the historical background 
of the particular case of the $(35_4)$ configuration. For further details, the reader is referred to~\cite{Bak, Klu, Lev, PS, Sal}; 
in a quite recent contribution, Conway and Ryba~\cite{CR} throw new light upon the \emph{Hexagrammum Mysticum}, i.e.\ 
the system  of 95 points and 95 lines consisting of the $(60_3)$ configuration of the Pascal lines and Kirkman points, and our
$(35_4)$ configuration.

We note as well that prior to Klug, a combinatorial configuration of type $(35_4)$ was presented in Brunel's work, as early as 1897--98.
Actually, Gr\"unbaum calls the attention to this work in his monograph~\cite{Gru09}. It turns out that this configuration is isomorphic 
with ours; moreover, as Gr\"unbaum remarks, it may be concluded that Brunel had found a geometric realization of this configuration
(cf.\ our Remark~\ref{ConclusionOnKlug} below).

Nevertheless, the first {\em known} geometric realization appears in fact in Klug's work~\cite{Klu}, Chapter 8. This point-line realization is visually unattractive. As we shall see later (in Sections~\ref{Constr} and~\ref{SectionCircles}), this situation can be improved.

\section[Relationship between DCD(n) and the Odd graph On]{Relationship between $\boldsymbol{\DCD(n)}$ and the\newline Odd graph $\boldsymbol{O_n}$} \label{Odd}

Recall that the {\em Kneser graph} $K(n, k)$ has as vertices the $k$-subsets of an $n$-element set, where two vertices 
are adjacent if the $k$-subsets are disjoint. The Kneser graph $K(2n -1,n - 1)$ is called an \emph{Odd graph} and is 
denoted by $O_n$~\cite{BH, God,GR}. By a simple comparison of this definition with that given in Section~\ref{Intro}, 
it is immediately clear that the configuration $\DCD(n)$ and the Odd graph $O_n$ are closely related to each other.

This relationship can be understood via the abstract $V$-construction, introduced and discussed in~\cite{GP}. Here we 
recall some definitions and results from that paper. 

Let $k\ge2$, $n\ge3$ be integers and let $G$ be a regular graph of valency $k$ on $n$ vertices. For a vertex $v$ of $G$, 
denote by $N(v)$ the set of vertices adjacent to $v$. Then take the family $S(G)$ of these vertex-neighbourhoods:
$$
S(G) = \{N(v) \,|\, v \in V (G) \},
$$
where $V(G)$ denotes the set of vertices of $G$.
The triple \hbox{$(V (G), S(G), \in)$} defines a combinatorial incidence structure, which we denote by $N(G)$. We call the graph 
$G$ {\em admissible\/} if no two of its vertices have a common neighbourhood (in \cite{PW} a graph with this property is called 
\emph{worthy}). For an admissible graph $G$, the incidence structure $N(G)$ is a combinatorial $(n_k)$ configuration. In this 
case we say that  $(n_k)$ is obtained from $G$ by {\em $V\!$-construction\/}. 

Recall that the {\em Levi graph\/} $L(C)$ of a configuration $C$ is a bipartite graph whose bipartition classes consist of the points 
and blocks of $C$, respectively, and two points in $L(C)$ are adjacent if and only if the corresponding point and block in $C$ are 
incident.  Many properties of configurations can be described by Levi graphs. For instance, a configuration is
flag-transitive if and only if its Levi graph is edge-transitive.
We have the following classical result~\cite{Cox50, Gru09, PS}.

\begin{lem} \label{Levi}
A configuration $C$ is uniquely determined by its Levi graph $L(C)$ with a given vertex coloring. 
\end{lem}

\noindent
If, in particular,  $C$ is obtained by $V$-construction, the following property is a useful tool in the study of such configurations~\cite{GP}.

\begin{thm} \label{LevKron}
Let $G$ be an admissible graph, and let $L$ be the Levi graph of the configuration $C$ obtained by $V\!$-construction from $G$. 
Then $L$ is the Kronecker cover of $G$. 
\end{thm}

\noindent
Here we recall that a graph $\tilde G$ is said to be the {\em Kronecker cover\/} (or {\em canonical double cover\/}) of 
the graph $G$ if there exists a $2:1$ surjective homomorphism $f : \tilde G \rightarrow G$ such that for every vertex $v$ 
of $\tilde G$ the set of edges incident with $v$ is mapped bijectively onto the set of edges incident with $f (v)$~\cite{GR, Por, IP}. 

We say that a configuration is {\em combinatorially self-polar\/} if there exists an automorphism of order two of its Levi graph 
interchanging the two parts of bipartition (see e.g.~\cite{PS}). The following result is also taken from~\cite{GP}.

\begin{thm} \label{Selfpolar}
A configuration obtained by $V\!$-construction from an admissible graph $G$ is combinatorially self-polar.
\end{thm}

Take now the Odd graph $O_n$, $n\ge2$, and consider a vertex $v$ of $O_n$. Since $v$ is an $(n-1)$-subset of the $(2n-1)$-element 
set, its neighbourhood $N(v)$ consists of $(n-1)$-subsets disjoint to $v$. Take the union of these latter sets, and denote it by $\bar N(v)$.
Clearly,  $\bar N(v)$ is complementary to $v$ in the $(2n-1)$-element set. It immediately follows that $O_n$ is admissible, hence the 
$V\!$-construction can be applied. It yields a configuration in which the block $\bar N(v)$ is incident with a point $v'$, such that $v'$
is adjacent to $v$ in the graph $O_n$. Hence the Levi graph of this configuration is a bipartite graph whose ``black" and ``white" vertices
are the $n$-subsets and $(n-1)$-subsets of the  $(2n-1)$-element set, respectively, where adjacency is given by containment. This 
graph is called the {\em revolving door graph\/}, or {\em middle-levels graph\/} (a subclass of the {\em bipartite Kneser graphs\/}; 
see e.g.~\cite{GP} and the references therein). On the other hand, it is directly seen that the Levi graph of the configuration $\DCD(n)$ 
is the same graph (consider the Cayley-Danzer arrangement from which it is derived). Hence, by Lemma~\ref{Levi}, we obtain the 
following main result of this section.
\begin{thm} \label{Vconst}
For $n\ge2$, the configuration $\DCD(n)$ is isomorphic to the combinatorial configuration obtained by $V\!$-construction 
from the Odd graph $O_n$. 
\end{thm}

\noindent
A consequence is that this combinatorial configuration has the same type $\left({\binom{2n-1}{n-1}}_n \right)$ as given in 
(\ref{DanzerType}) in Section~\ref{Intro}. The configuration $\DCD(2)$ corresponds to the trilateral $(3_2)$, the case $n=3$ 
gives rise to the Desargues con\-fig\-u\-ration $\DCD(3)$ of type $(10_3)$, while the particular case $n=4$ corresponds to 
Danzer's configuration $(35_4)$.

\begin{cor}
Danzer's configuration $(35_4)$ is isomorphic to the combinatorial configuration obtained by $V\!$-construction from the Odd graph $O_4$. 
\end{cor}

\begin{figure}[h!] 
  \begin{center} 
  \includegraphics[width=0.85\textwidth]{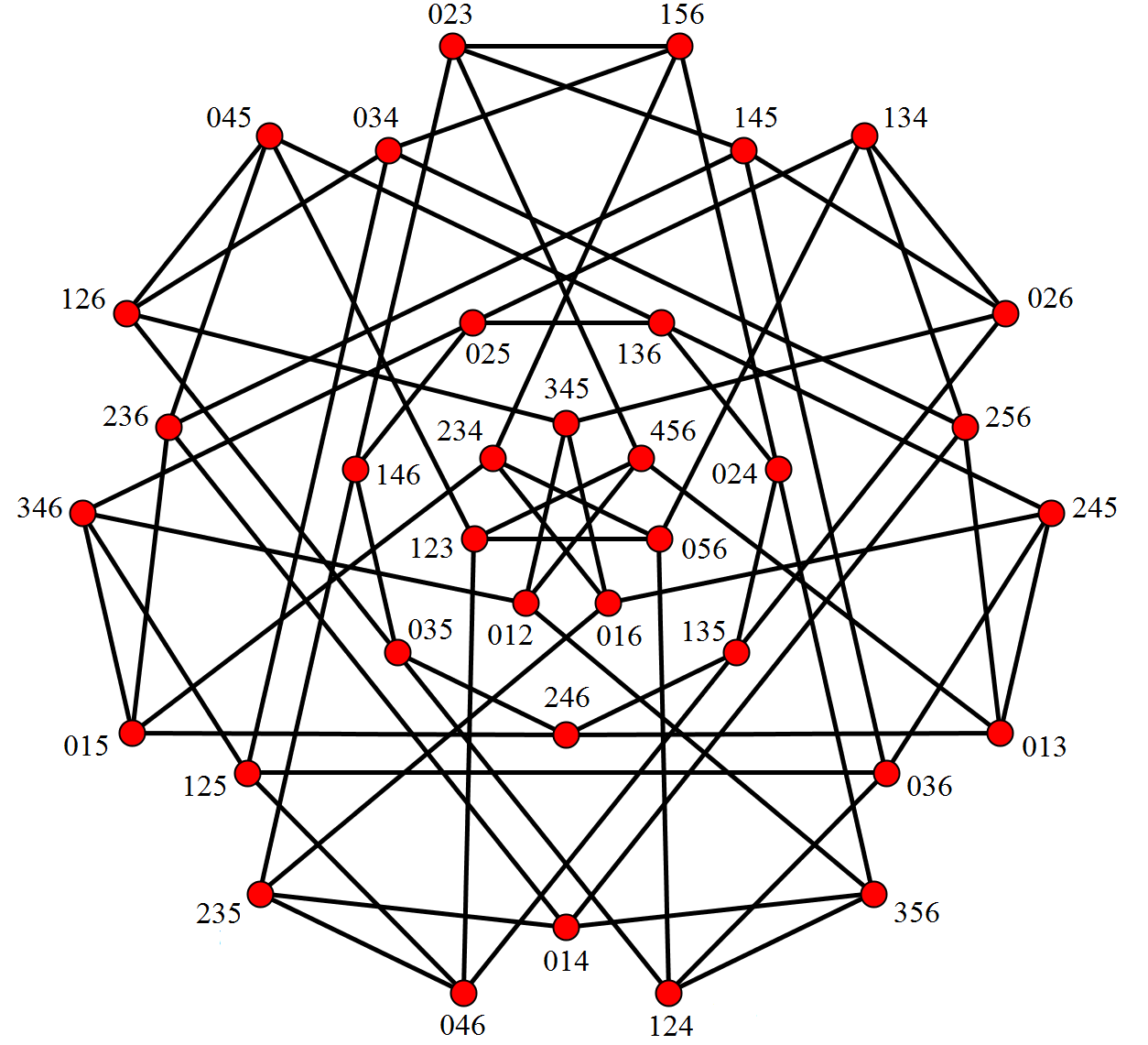} 
  \caption{The Odd graph $O_4$, realized with a 7-fold dihedral symmetry.}
  \label{Odd4} 
  \end{center} 
\end{figure}

\noindent
A representation of the graph $O_4$ is depicted in Figure~\ref{Odd4}.

\section[Decomposition of DCD(n)]{Decomposition of $\boldsymbol{\DCD(n)}$} \label{Decomp}

The following definition is taken from~\cite{GG14}.

\begin{definition}\label{IncSum}
By the {\it incidence sum\/} of configurations $\mathcal C_1$ and $\mathcal C_2$ we mean the configuration $\mathcal C$,
which is the disjoint union of $\mathcal C_1$ and $\mathcal C_2$, together with a specified set 
$I \subseteq \mathcal P_1 \times \mathcal L_2 \cup  \mathcal P_2 \times \mathcal L_1$ of incident point-line pairs, where
$\mathcal P_i$ denotes the point set and $\mathcal L_i$ denotes the line set of $\mathcal C_i$, for $i = 1,2$. 
We denote it by $\mathcal C_1\oplus_I \mathcal C_2$.
\end{definition}

\begin{figure}[h!]
  \begin{center}  
   \includegraphics[width=0.85\textwidth]{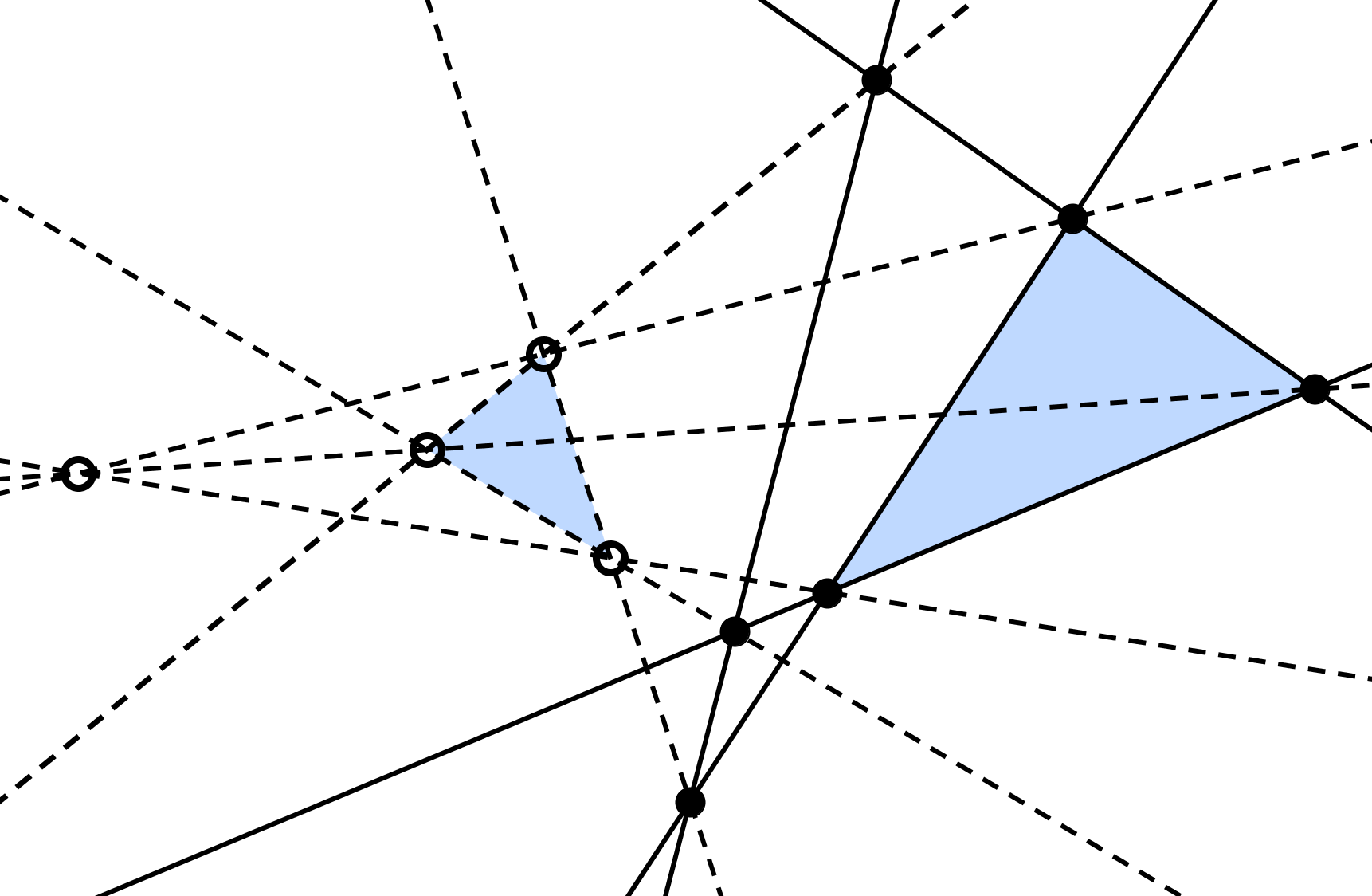}
    \caption{The Desargues configuration $(10_3)$ is the incidence sum of a complete quadrangle $(4_3, 6_2)$ 
                   and a complete quadrilateral $(6_2, 4_3)$. (The former denoted by white points and dashed lines, 
                   while the latter by black points and solid lines.)}
    \label{DesDec}
    \end{center}
\end{figure}

\noindent
If the set $I$ is clear from the context, it can be omitted from the operation symbol.  Figure~\ref{DesDec} shows the example that the 
Desargues configuration can be considered as the incidence sum of a complete quadrangle and a complete quadrilateral (or, in other 
words, the incidence sum of the \emph {Pasch configuration} and its dual). The set of new incidences consists of 6 point-line pairs 
formed by the lines of the complete quadrangle and by the points of the complete quadrilateral.  

We show that this is a particular (actually, the second smallest) case of the following general relationship.

\begin{thm} \label{ThmDecomp}
For all $n\ge3$, the configuration $\DCD(n)$ is the incidence sum of the form $\mathcal C_1\oplus_I \mathcal C_2$ such that
\begin{enumerate}[(1)]
\item $\mathcal C_1$ is a configuration of type
\begin{equation}\label{C1type}
\left({\binom{2n-2}{n-2}\!}_n, {\binom{2n-2}{n-1}\!}_{n-1}\,\right );
\end{equation}
\item $\mathcal C_2$ is a configuration of type
\begin{equation}\label{C2type}
\left({\binom{2n-2}{n-1}\!}_{n-1}, {\binom{2n-2}{n-2}\!}_n \,\right);
\end{equation}
\item the set I of new incidences consists of $\,\binom{2n-2}{n-1}$ point-line pairs whose points belong to $\mathcal C_2$ 
and whose lines belong to $\mathcal C_1$;
\item $\mathcal C_1$ and $\mathcal C_2$ are dual to each other; 
\item $\mathcal C_1$ and $\mathcal C_2$ are flag-transitive configurations. 
\end{enumerate}
\end{thm} 

\noindent
\emph{Proof.}
Consider first the $\binom{2n-1}{n-1}$ points which originate from the $\binom{2n-1}{n-1}=\binom{2n-1}{n}$ $n$-tuples 
of hyperplanes determining $\DCD(n)$, as described in Section~\ref{Intro}. Using Pascal's rule:
\begin{equation}\label{eq:Pascal}
\binom{2n-1}{n} =  \binom{2n-2}{n} +  \binom{2n-2}{n-1},
\end{equation}
we see that this set of $n$-tuples decomposes into the disjoint union of a set  consisting of $\binom{2n-2}{n}=\binom{2n-2}{n-2}$ 
$n$-tuples and a set consisting of $\binom{2n-2}{n-1}$ $n$-tuples. Denote these sets by $P_1$ and $P_2$, respectively. For an 
arbitrary but fixed hyperplane $H$, they can be chosen so that the $n$-tuples in $P_2$ contain, but the $n$-tuples in $P_1$ do not  
contain  $H$. Let $L_1$ be the set of $(n-1)$-tuples obtained from the $n$-tuples of $P_2$ by omitting $H$ from them. Let 
\begin{equation}\label{phi}
\varphi: L_1 \rightarrow P_2
\end{equation} 
be the bijection defined by this procedure in a natural way. Furthermore, let  $L_2$ be the set complementary to $L_1$ in the set of  
$(n-1)$-tuples of hyperplanes determining the lines of $\DCD(n)$. 

Thus, $P_1$ and  $L_1$ consist of $\binom{2n-2}{n-2}$ $n$-tuples and $\binom{2n-2}{n-1}$ $(n-1)$-tuples, respectively, such that 
they do not contain $H$. On the other hand,  $P_2$ and  $L_2$ consist of $\binom{2n-2}{n-1}$ $n$-tuples and $\binom{2n-2}{n-2}$
$(n-1)$-tuples, respectively, each containing $H$.

Hence, $P_1$ and $L_1$ yield a set of points and a set of lines, respectively, with cardinalities corresponding to type (\ref{C1type}); 
the same is valid in the case of $P_2$ and $L_2$, regarding type (\ref{C2type}). Besides, from each of the $n$-tuples in $P_1$ one 
can obviously omit precisely $n$ distinct hyperplanes so as to obtain an $(n-1)$-tuple; thus, each point determined by these $n$-tuples 
is incident with precisely $n$ lines determined by the $(n-1)$-tuples obtained in this way. Note that, since the $n$-tuples do not contain 
$H$, neither do the $(n-1)$-tuples; hence the set of all these latter is just equal to $L_1$. Furthermore, each of these $(n-1)$-tuples can 
be completed by a hyperplane  other than $H$ in precisely $(n-1)$ distinct ways. Thus we obtain that the points and lines yielded by the 
sets $P_1$ and $L_1$, respectively, form a configuration $\mathcal C_1$ of the desired type (\ref{C1type}). By similar reasoning, one 
sees that a configuration $\mathcal C_2$ of type (\ref{C2type}) can be obtained from the sets $P_2$ and $L_2$.

The pairs of points and lines determined by the $n$-tuples and $(n-1)$-tuples, respectively, which correspond to each other under the 
bijection $\varphi$ in (\ref{phi}) form the set $I$ of (\ref{C2type}).

Now we define the map $\delta: \mathcal C_1 \rightarrow \mathcal C_2$, as follows.  For any point $p$ in $\mathcal C_1$ determined by 
an $n$-tuple of hyperplanes, take the line $l$ in $\mathcal C_2$ determined by the complementary $(n-1)$-tuple, and set $\delta (p) = l$.
Likewise, for any line $m$ in $\mathcal C_1$ determined by an $(n -1)$-tuple of hyperplanes, take the point $q$ in $\mathcal C_2$ 
determined by the complementary $n$-tuple, and set $\delta (m) = q$. De Morgan's laws imply that $\delta$ is a duality map. The 
symmetry of the construction which does not depend on the order of the elements in $n$-tuples enables us to conclude that both 
configurations are flag-transitive.
\hfill $\blacksquare$
\medskip

Here we emphasize again the fact that has already been touched following Theorem~\ref{TripleDesDualThm}.

\begin{remark} \label{NoExtension}
The duality map $\delta$ just defined cannot be extended to include, in the particular case of $n=4$, the Pascal lines and Kirkman points
as well. For, although they form a configuration $(60_3)$~\cite{Lev, PS}, this configuration is not combinatorially self-polar~\cite{Klu, CR}.
\end{remark}

\begin{cor}
The configuration $\DCD(n)$ is self-polar.
\end{cor}

\noindent
\emph{Proof.} The statement already follows from Theorems~\ref{Vconst} and~\ref{Selfpolar}. A self-polarity map $\pi$ can now be 
realized as follows. Consider a decomposition $\DCD(n) = \mathcal C_1\oplus_I \mathcal C_2$ as in Theorem~\ref{ThmDecomp}, and 
take the duality map $\delta: \mathcal C_1 \rightarrow \mathcal C_2$. Then
$$
\pi(x) = \left\{
\begin{array}{ll}
  \delta(x) & \mbox {if $x$ is a point or line in $\mathcal C_1$,}\\
  \delta^{-1}(x) & \mbox {if $x$ is a point or line in $\mathcal C_2$.}
\end{array}
\right. 
$$ 
\vskip -10pt
\hfill $\blacksquare$

\bigskip

In the particular case of $n=3$, the decomposition theorem~\ref{ThmDecomp} leads to the following rephrasing of Desargues' theorem.

\begin{thm}
The following statement is equivalent to Desargues' theorem. Given a complete quadrangle $Q$, one can choose six points, one on each 
line of $Q$, such that these points determine a complete quadrilateral.
\end{thm}

\noindent
\emph{Proof.} First we show that Desargues' theorem implies the statement. Consider the Desargues configuration, in which a complete 
quadrangle $Q$ is chosen, as it is depicted in Figure \ref{DesDec}. Note that the centre of perspectivity is the white point not incident 
with the shaded triangle of $Q$. Let $Q'$ denote the desired quadrilateral. For finding the vertices of $Q'$, first choose the three black 
points which are vertices of the other shaded triangle. Then we choose three more points, such that each is formed as the intersection of 
the lines of the two triangles corresponding to each other under central perspective. But these latter three points are collinear by Desargues' 
theorem, so the six points altogether form indeed a complete quadrilateral $Q'$.

The converse implication can simply be obtained by taking the complete quadrangle and the complete quadrilateral as in Figure \ref{DesDec}, 
and choosing the shaded triangles as  corresponding to each other under perspectivity. 
$\blacksquare$
\bigskip

\noindent 
Note that likewise to Desargues' theorem, the dual of our equivalence statement also holds.

In what follows we apply the decomposition theorem~\ref{ThmDecomp} to Danzer's configuration $(35_4)$.

\begin{lem}\label{SP_Lemma}
The Steiner-Pl\"ucker configuration $(20_3, 15_4)$ in Theorem~\ref{TripleDesThm} can be obtained in the following way. Take 
seven hyperplanes, $H_0,\dots,H_6$ in general position in the 4-dimensional projective space. Take all the quadruples and triples 
such that each contain, say, $H_0$. The intersections of the hyper\-planes in these subsets yield 20 points and 15 lines, respectively, 
which, by suitable projection into the projective plane, form the desired configuration.
\end{lem}

\begin{figure}[h!]
  \begin{center}  
   \includegraphics[width=.95\textwidth]{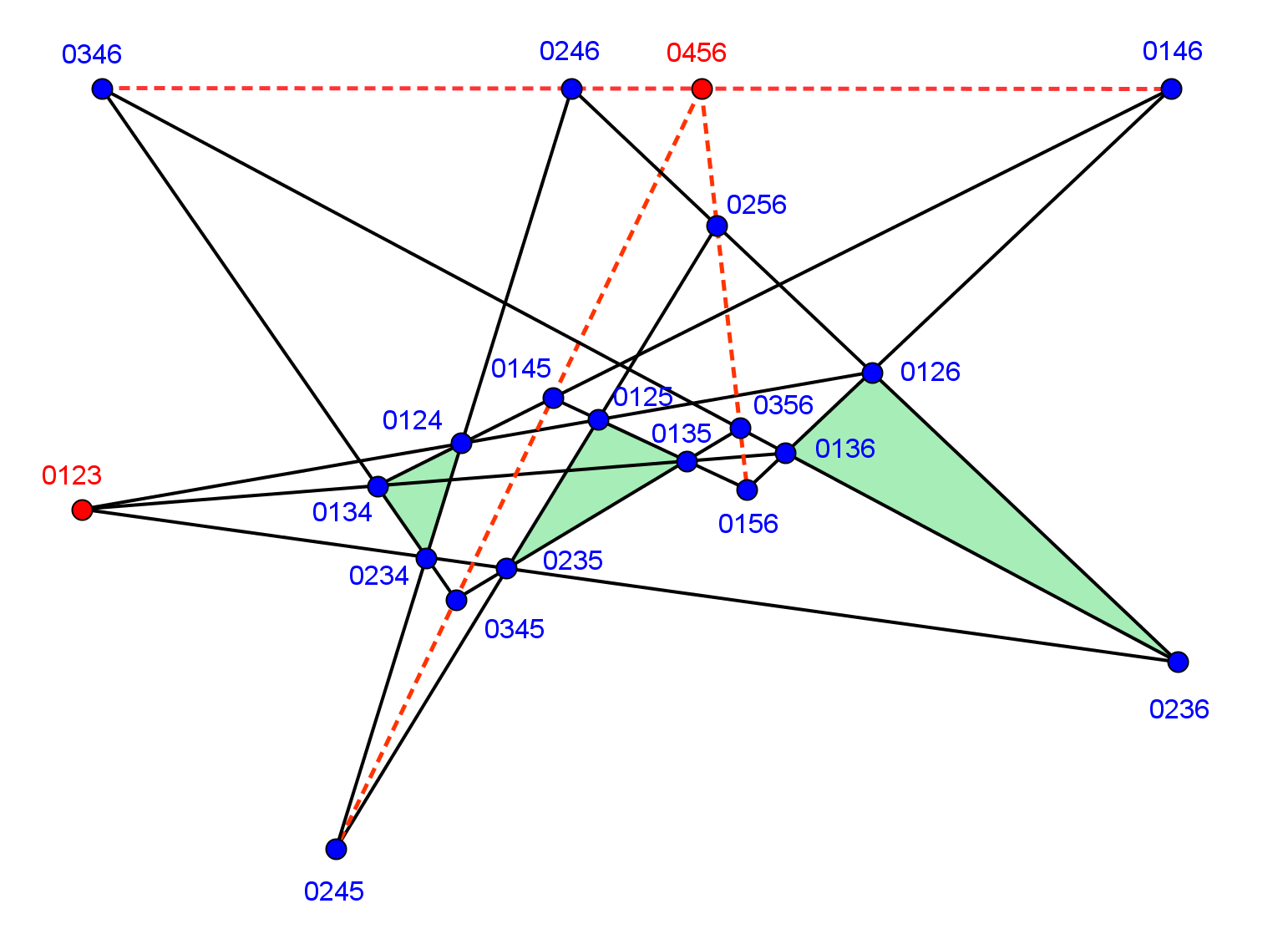}
    \caption{The Steiner-Pl\"ucker configuration $(20_3, 15_4)$. The points are labelled in accordance with Lemma~\ref{SP_Lemma}. 
                   The shaded triangles are perspective from the same point (labelled by $0123$), as in Theorem~\ref{TripleDesThm}; 
                   the point of concurrency of the three axes of perspectivity is labelled by $0456$.}
    \label{SP_Labels}
    \end{center}
\end{figure}

\noindent
\emph{Proof.} See the labelling in Figure~\ref{SP_Labels}, which refers to the quadruples of hyperplanes determining the points. 
The intersection of the labels of points incident with the same line corresponds to the triple of hyperplanes determining that line.
\hfill $\blacksquare$
\vskip 10pt

\noindent
Thus, together with that arising from Theorem~\ref{TripleDesThm}, one can see here three different ways of deriving the 
Steiner-Pl\"ucker configuration. We remark that Adler gives a fourth one, see~\cite{Adl}. Independently, this configuration 
is also presented by Servatius and Servatius~\cite{SS} under the name of \emph{generalized Reye configuration}.
It is possible to verify that our Figure \ref{SP_Labels} depicts the same configuration as that in Figure 5 of \cite{SS}; 
see also Figure 6.50 in~\cite{PS}. 

We also note the following interesting realization of this configuration. Szilassi observed~\cite{Szi} that in addition 
to the original three triangles there is another triple of triangles which seems to be in some definite relationship with the former. 
Using the labelling in our Figure~\ref{SP_Labels}, it readily follows that there is a combinatorial automorphism $\varphi$ 
of the configuration which, in terms of the defining hyperplanes in Lemma~\ref{SP_Lemma}, can be given in the form
\begin{align*}
1 & \longleftrightarrow 6 \\
2 & \longleftrightarrow 5 \\
3 & \longleftrightarrow 4. \\
\end{align*}
\vskip -18pt

\noindent
In general, $\varphi$ cannot be realized geometrically (in the most general case, by a projective collineation of  period two, 
i.e.\ by a \emph{harmonic homology}~\cite{Cox74}); in other words, given a particular realization of the Steiner-Pl\"ucker 
configuration, there is no geometric transformation carrying it into itself that would induce $\varphi$. However, we found that in some special cases
there exist even centrally symmetric realizations. 
\begin{lem}\label{CS_Lemma}
The Cayley-Salmon configuration $(15_4, 20_3)$ in Theorem~\ref{TripleDesDualThm} can be obtained in the following way. Take six
hyperplanes, $H_1,\dots,H_6$, in general position in the 4-dimensional projective space.  They meet by fours and by threes in 15 points 
and 20 lines, respectively. The configuration formed by these points and lines is then to be suitably projected into the projective plane. 
\end{lem}

\noindent
\emph{Proof.} Since the Steiner-Pl\"ucker configuration and the Cayley-Salmon configuration are dual to each other, the statement is a 
straightforward consequence of Theorem~\ref{ThmDecomp} and Lemma~\ref{SP_Lemma}. 
\hfill $\blacksquare$
\bigskip

\begin{figure}[h!]
  \begin{center}  
   \includegraphics[width=1\textwidth]{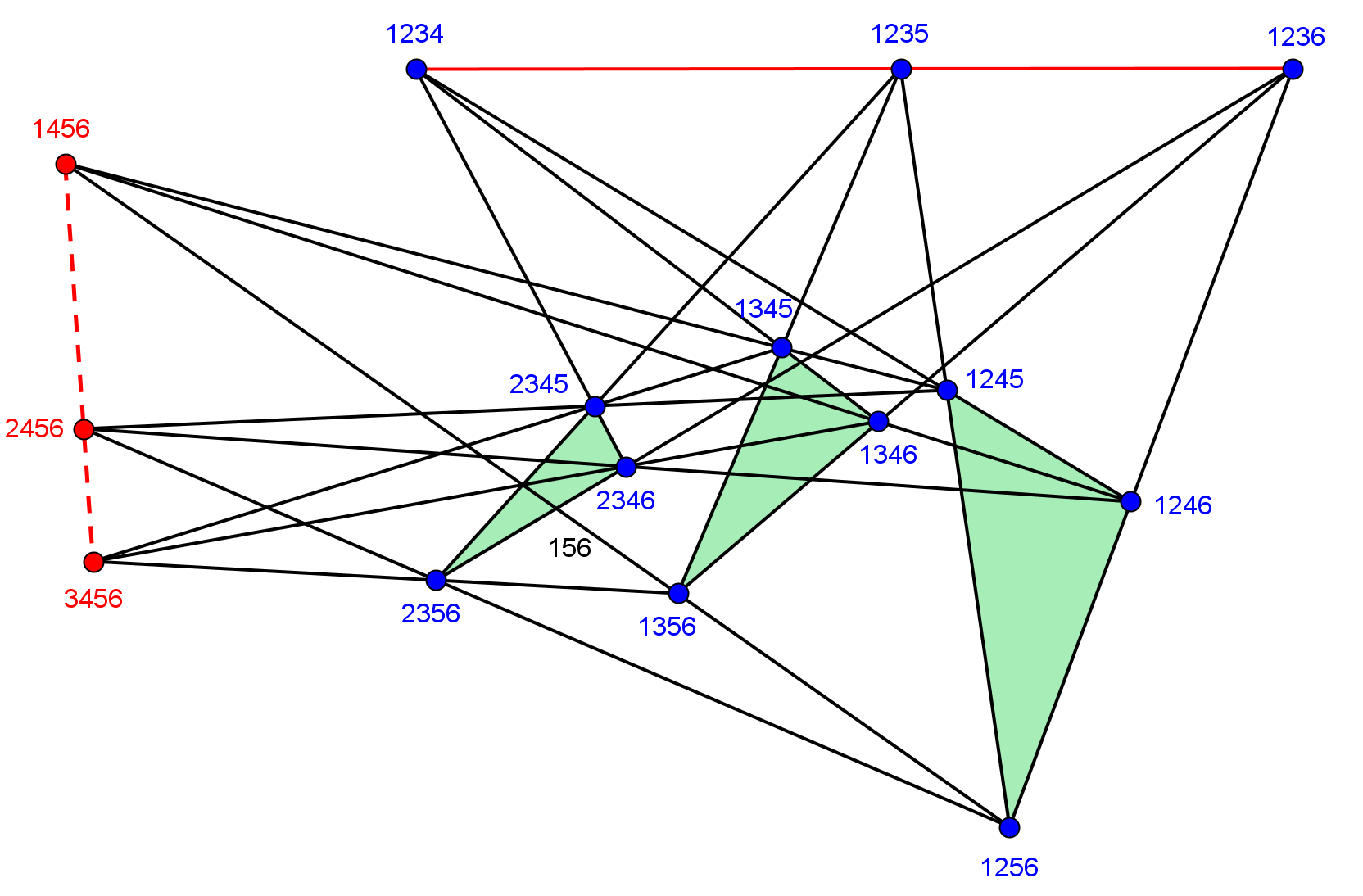}
    \caption{The Cayley-Salmon configuration $(15_4,20_3)$. The points are labelled in accordance with Lemma~\ref{CS_Lemma}; 
                   they also indicate the duality with the Steiner-Pl\"ucker configuration depicted in Figure~\ref{SP_Labels}, as it is described 
                   in the proof of Theorem~\ref{ThmDecomp}.}
    \label{CS_Labels}
    \end{center}
\end{figure}

\noindent
Note that similarly to the proof of the previous lemma, this latter can also be inferred from the labelling of Figure~\ref{CS_Labels},
which depicts the Cayley-Salmon configuration.

\begin{thm} \label{DanzerDecomp}
Danzer's configuration $(35_4)$ can be decomposed into the incidence sum of the Steiner-Pl\"ucker configuration $(20_3, 15_4)$ and 
the Cayley-Salmon configuration $(15_4, 20_3)$.
\end{thm}

\noindent
\emph{Proof.} The statement follows from Theorem~\ref{ThmDecomp}, by Lemma~\ref{SP_Lemma} and Lemma \ref{CS_Lemma}. 
\hfill $\blacksquare$
\begin{remark} \label{ConclusionOnKlug}
A comparison of Theorems~\ref{KlugObservation} and ~\ref{DanzerDecomp} shows that Klug knew the 
configuration $(35_4)$ what we call Danzer's configuration (to be precise, we traced it back to this point; cf.\ Remark (vi) 
of Gr\"unbaum in his ``Musings"~\cite{Gru08}).
\end{remark}

\begin{remark} Poto\v{c}nik proved \cite{P13} that there is exactly one biregular graph with 20 vertices of valence 3 and 15 vertices 
of valence 4 that is edge-transitive of girth 6. In his Table 2 it has $ID = \{35,2\}$. One may conclude that this graph is the Levi graph 
of both the Cayley-Salmon and the Steiner-Pl\"{u}cker configuration. From \cite{PSV} it follows that there are only two flag-transitive 
self-polar combinatorial $(35_4)$ configurations. From \cite{census} we may deduce that the other one is a cyclic configuration 
with the symbol $\{0,1,8,14\}$. In particular, this means its Levi graph is a cyclic Haar graph, see \cite{HMP}, and is denoted as 
{\tt C4[70,3]} in \cite{census}. The \emph{Danzer graph}, i.e.\ the Levi graph of Danzer's configuration is denoted as {\tt C4[70,4]}.
It is the Kronecker cover over the Odd graph $O_4$ (cf.\ Theorem 3.2 and Corollary 3.5).

We refer the reader to an interesting paper \cite{rotary} in which further information about
the nature of both configurations can be read off the Table 4 for $v = 35$. In that table 
{\tt C4[70,3]} is immediately followed by {\tt C4[70,4]}.
  
\end{remark}

\section {Another subconfiguration of Danzer's configuration} \label{Coxeter}

Here we show that besides those given in Theorem~\ref{DanzerDecomp} there is a further interesting subconfiguration of Danzer's 
$(35_4)$. It played some role in our attempts to find a realization of Danzer's configuration with seven-fold rotational symmetry.

\begin{figure}[h] 
  \centering\mbox{\hskip -3pt
    \subfigure[]{
    \includegraphics[width=0.475\textwidth]{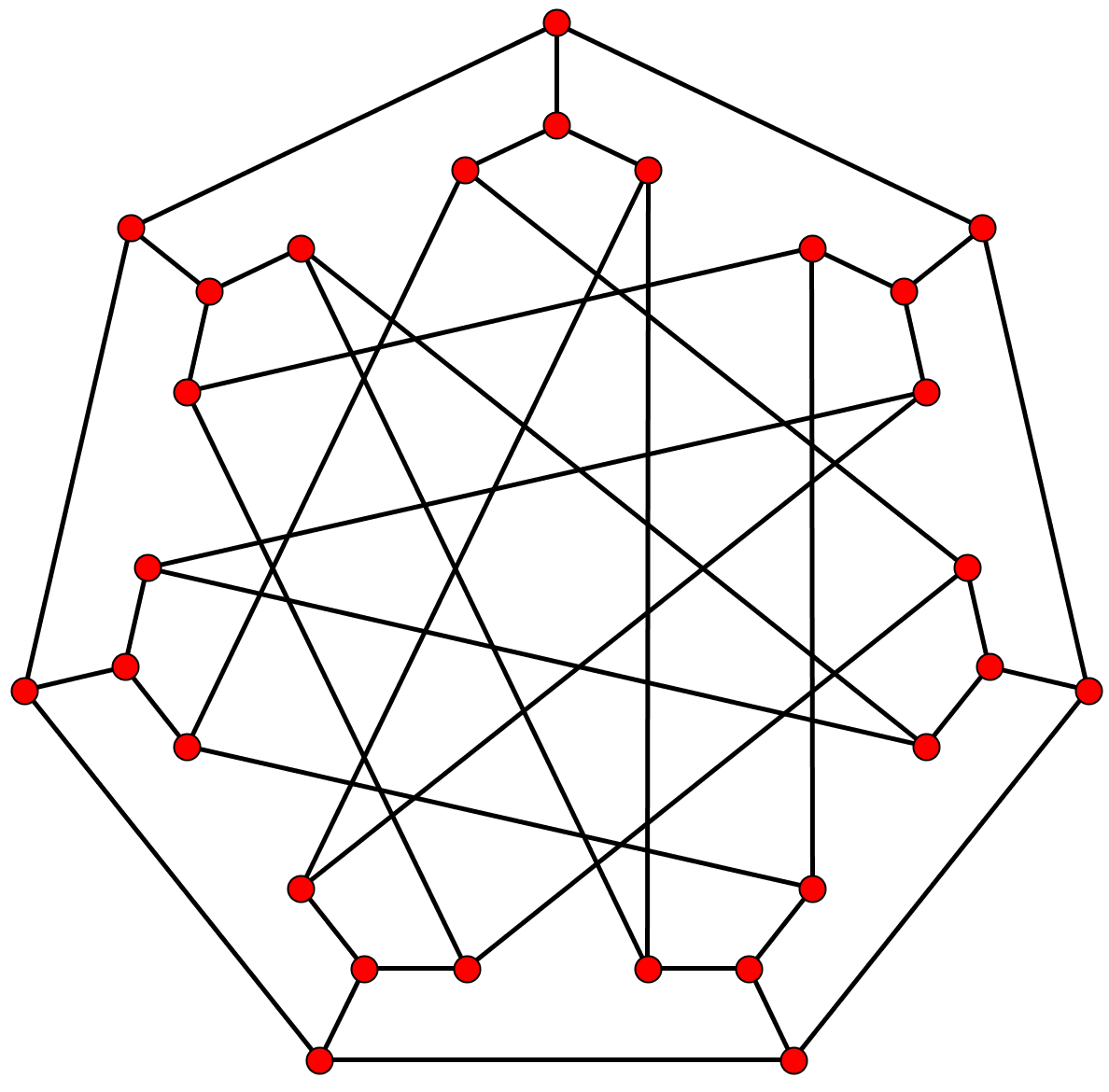}} \quad
    \subfigure[]{
    \includegraphics[width=0.475\textwidth]{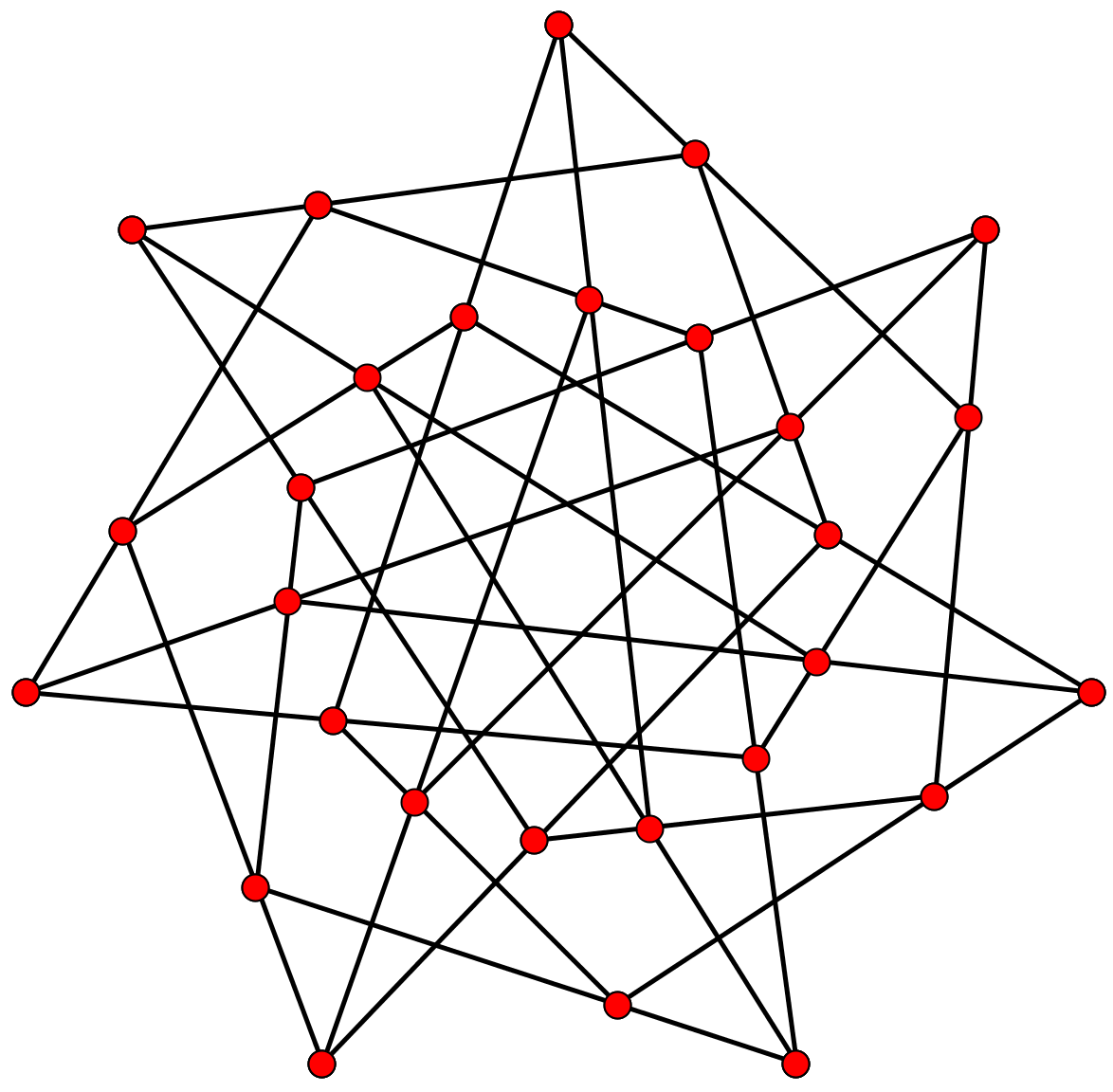}}}
    \caption{The Coxeter graph (a), and the corresponding point-line Coxeter  $(28_3)$-configuration  (b), both with
                   7-fold rotational symmetry.}                
    \label{CoxGraphConf}
\end{figure}

Consider the Coxeter graph, which is a 28-vertex cubic graph. A particularly nice drawing of it, due to Randi\' c~\cite{BM}, is shown 
in Figure~\ref{CoxGraphConf}a. The vertices can be identifed with the 28 antiflags, i.e. non-incident point-line pairs $(p,l)$, of the 
Fano plane $\mathcal F$. Two vertices $(p,l)$ and $(p',l')$ are adjacent if and only if $p\ne p'$, $l \ne l'$, and the pairs $(p,l')$, $(p',l)$ 
are also non-incident~\cite{God}. This is equivalent to saying that the vertices are labelled with triplets of points in $\mathcal F$ 
which do not lie on the same line, and adjacent vertices are labelled with disjoint triplets. It follows that the Coxeter graph is a 
subgraph of the Kneser graph $K(7, 3)$~\cite{GR}.

A point-line realization of the combinatorial configuration $(28_3)$ obtained from the Coxeter graph $\Gamma$ by $V\!$-construction 
is shown in Figure~\ref{CoxGraphConf}b. We denote it by $N(\Gamma)$ and name it the \emph{Coxeter $(28_3)$-configuration}.

\begin{remark} 
Sometimes a configuration of type $(12_3)$ that has the Nauru graph, i.e. the arc-transitive generalized Petersen graph on 24 vertices, 
as its Levi graph, is called the Coxeter configuration, instead of Nauru configuration. If we want to make distinction between the two, 
we propose that the latter be called \emph{Coxeter $(12_3)$-configuration}.
\end{remark} 

On the basis of the close structural relationship between the Coxeter graph and the Odd graph $O_4$ we expected that the rotational 
realization of $N(\Gamma)$ could be extended to a heptagonally symmetric realization of Danzer's configuration. However, all our 
attempts to find such a realization failed. We note that such an attempt is not quite new. In \cite{Gru08} Gr\"{u}nbaum admits that 
he unsuccessfully tried to find a symmetric drawing of Danzer's configuration. Therefore, we suggest the following conjecture.

\begin{conj} \label{NonReal}
There is no realization of Danzer's $(35_4)$ point-line configuration with five- or seven-fold rotational symmetry.
\end{conj}

\noindent
We remark that this is in accordance with a similar property of Desargues' $(10_3)$ point-line configuration, which can be 
proved~\cite{GP}. One may also ask whether such a property of non-realizability of $\DCD(n)$ with $n$-fold symmetry  is
valid for all $n$. We did not investigate this problem in such generality; even verification of the conjecture for $n=4$ needs 
further work (possibly with detailed calculations), which is beyond the scope of the present paper.

\section[Constructing point-line realization of Danzer's configuration]{Constructing point-line realization of\\ Danzer's configuration} \label{Constr}

\noindent
Our first construction is based on Theorem~\ref{DanzerDecomp}. Start from a concrete realization $\mathcal C_1$ of the 
Steiner-Pl\"ucker configuration $(20_3, 15_4)$, e.g.\ from that depicted in Figure \ref{SP_Labels}. Choose the complete 
quadrangle determined  by the vertices $0123, 0124, 0134, 0234$. Take a complete quadrilateral determined by four new 
lines such that they are labelled by $456, 356, 256, 156$, and they are incident with the vertices $0456, 0356, 0256, 0156$, 
respectively. By a basic theorem of projective geometry~\cite{Cox74}, there is a unique projective correlation sending the 
vertices of the quadrangle into the lines of the corresponding quadrilateral. In our case this correlation is just the duality 
map $\delta$ defined in the proof of Theorem~\ref{ThmDecomp}, and determined by the labels chosen here. Accordingly, 
the complete quadrilateral can be extended to a realization $\mathcal C_2$ of the Cayley-Salmon configuration 
$(15_4, 20_3)$. The incidence sum $\mathcal C_1\!\oplus \mathcal C_2$ yields the Danzer's configuration.

It is directly seen that one cannot expect that this construction  would yield a visually attractive realization at all. Hence, 
we prefer the following one.

In~\cite{BP} a concept of polycyclic configurations was introduced. A configuration is {\em polycyclic\/} if there exists a non-trivial 
cyclic automorphism $\alpha$ which is semi-regular (all orbits of points and lines are of the same size). If $\alpha$ has order $k$, 
then such configuration is called $k$-cyclic. If a realization of a $k$-cyclic configuration in the Euclidean plane realizes $\alpha$ 
as a rotation for angle $2\pi/k$, then such realization (if exists) is called {\em rotational}. 

A theory was developed in~\cite{BP} 
to analyze polycyclic configurations via voltage graphs.

As usual, let $E(G)$ denote the set of edges of $G$. By $A(G)$ we denote the set of {\em arcs} of  $G$, i.e.\ each edge $e = uv$ determines two opposite arcs $a= (u,v)$ and $ a^* =(v,u)$. 
If $a$ is an arc, then $a^*$ represents the opposite arc. Clearly $a \mapsto a^*$ is an involution on $A(G)$, and the edges of $G$ may be viewed as orbits of this involution.
A {\em voltage graph\/} is a triple $(G, \Gamma, \alpha)$ where $G$ is a directed graph (digraph), $\Gamma$ is a group, and $\alpha : A(G) \to \Gamma$ is a map assigning a {\em voltage\/} to each of the arcs such that $\alpha (a^*) = \alpha^{-1}(a)$, for any arc $a$.
 The {\em covering graph\/} over a voltage graph $(G, \Gamma, \alpha)$ is the graph $\tilde{G}$ with the vertex set $V(\tilde{G}) = V(G) \times \Gamma$ and the edge set
 $E(\tilde{G}) = \{(u,g)(v,g + \alpha(a)) : (u,v) = a \in A(G), g \in \Gamma\}$. Note that each edge of $\tilde G$ is defined twice, but the definition is consistent.

 Also, $\Gamma$ will always be a cyclic group $\mathbb{Z}_k$ and for the clarity of figures we will omit labels and directions of the edges with voltage $0$. For more details on voltage graphs see, e.g.,~\cite{GT}. 

It is easy to see~\cite{BP} that the Levi graph of a $k$-cyclic configuration is isomorphic to a $\mathbb{Z}_k$ covering graph over some (bipartite) voltage graph (called the {\em reduced Levi graph\/} in~\cite{Gru09}). Conversely, a covering graph of a bipartite voltage graph with group $\mathbb{Z}_k$ is a Levi graph of some $k$-cyclic configuration whenever no parallel pair of edges, no two adjacent pairs of edges and no $4$-cycle in the voltage graph lift to a cycle of length $2$ or $4$.

\begin{figure}[h] 
  \begin{center}\hskip 4pt
  \subfigure[]{\hskip -4pt 
  \includegraphics[width=0.35\textwidth]{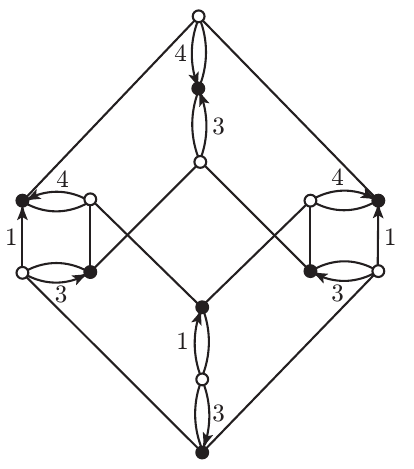}} \hskip 25pt
 \subfigure[]{\hskip -4pt 
   \includegraphics[width=0.35\textwidth]{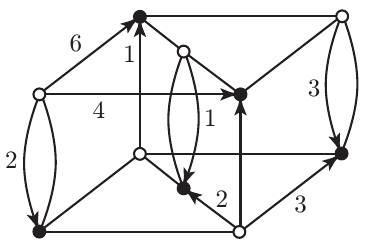}}
  \end{center}
  \caption{Voltage graphs for polycyclic realization of Danzer's configuration; $5$-cyclic (a), and $7$-cyclic (b).}
  \label{fig:voltage}  
\end{figure}
\begin{figure}[h!]
\centering
\includegraphics[width=0.85
\textwidth]{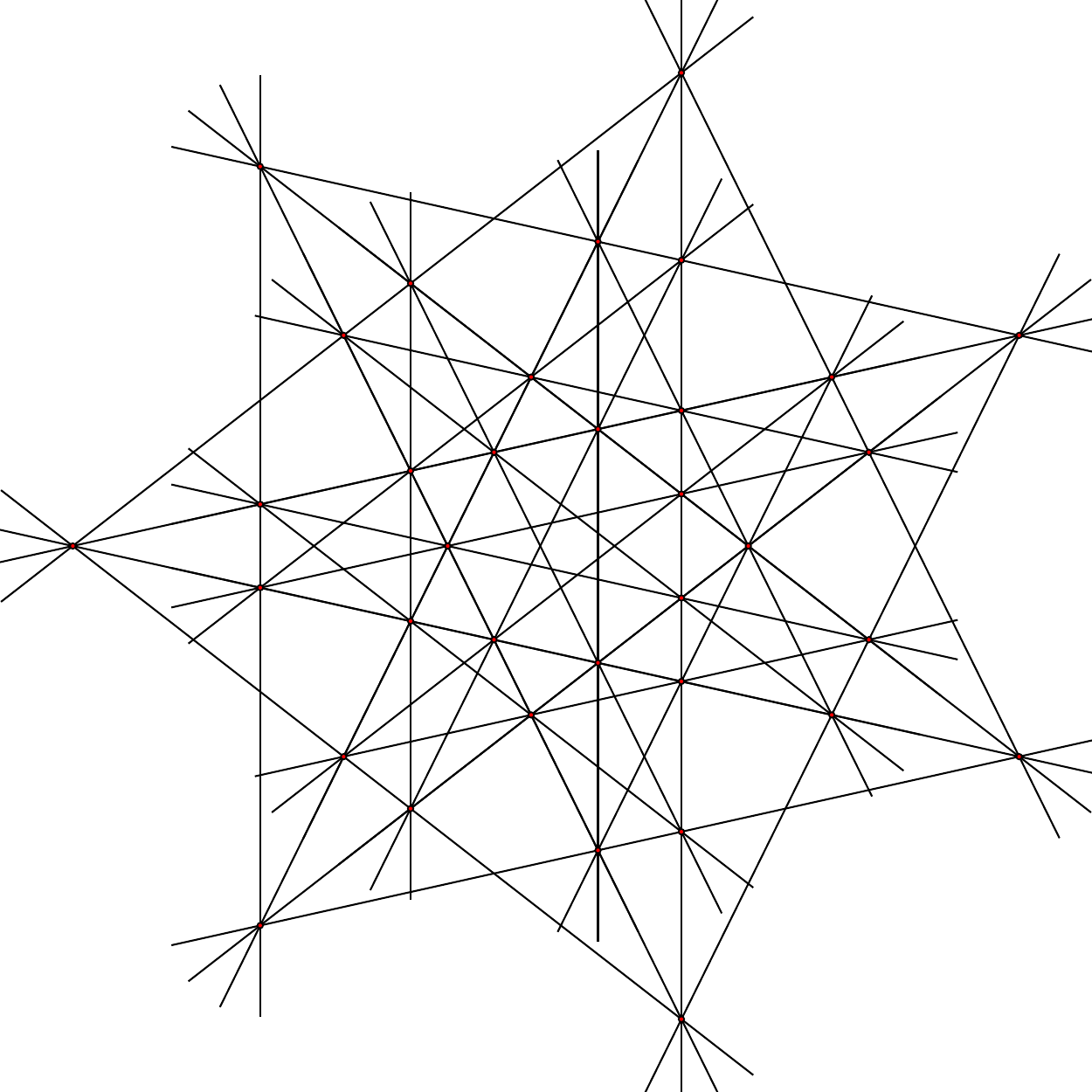}
\caption{``Realization" of Danzer's configuration with 7-fold dihedral symmetry, where lines from two different orbits coincide. }
\label{fig:coincide}
\end{figure}
\begin{figure}[h!]
\centering
\includegraphics[width=0.85\textwidth]{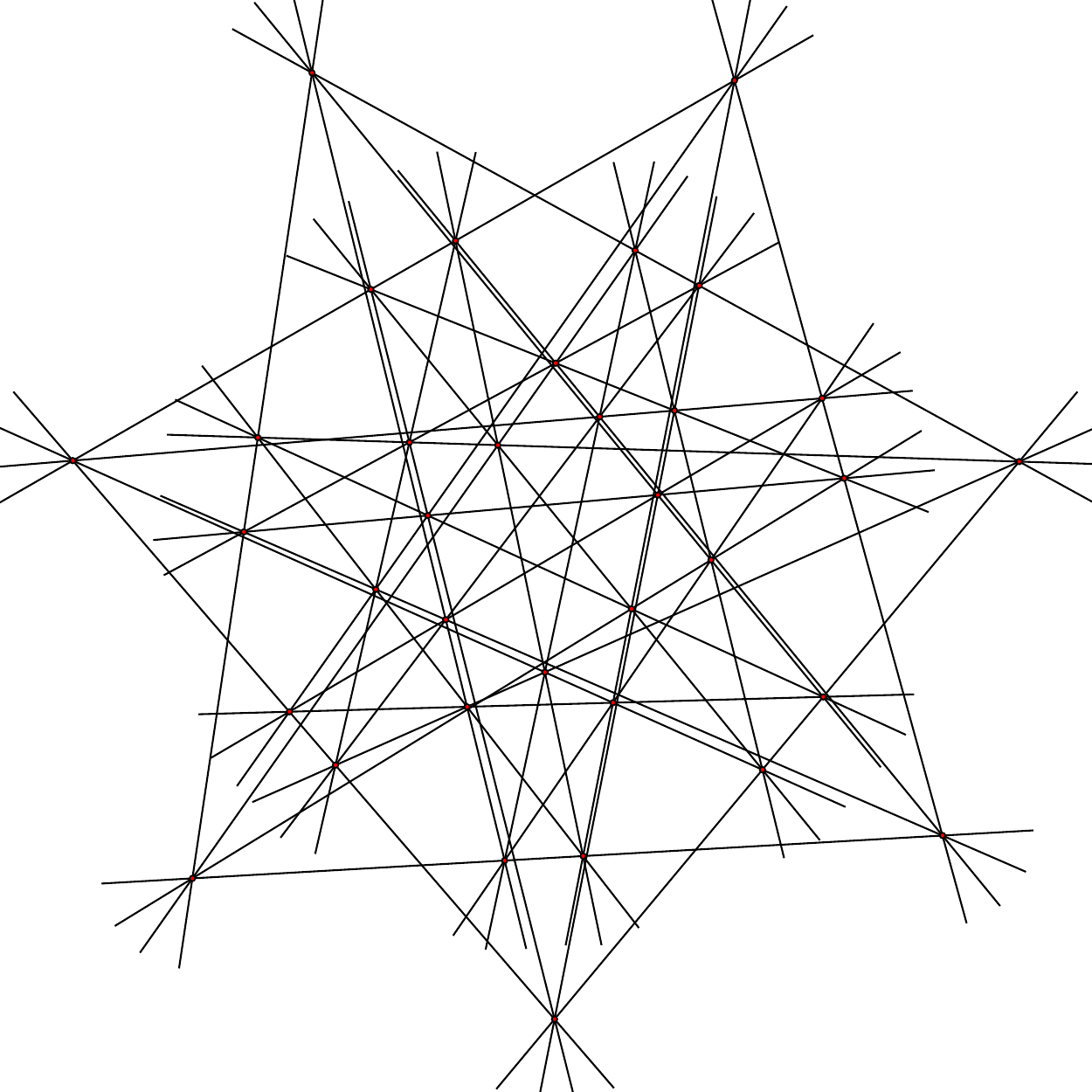}
\caption{Realization of Danzer's configuration resembling the $7$-fold rotational symmetry.}
\label{fig:perturb}
\end{figure}

Danzer's configuration is $5$-cyclic and $7$-cyclic.  The automorphism group of its Levi graph is $S_7 \times \mathbb Z_2$, 
the direct product of the symmetric group of degree 7 with the group of order two; hence its order is 10\,080. The Levi graph is, 
for example, a $\mathbb{Z}_5$ covering graph over the voltage graph in Figure~\ref{fig:voltage}a, and a $\mathbb{Z}_7$ 
covering graph over the voltage graph in Figure~\ref{fig:voltage}b. To obtain rotational realization it is necessary to find 
solutions of polynomial equations which depend directly on the voltage graph structure. Using {\sc Mathematica} we have 
analyzed all possible voltage graphs for $5$- and $7$-cyclic structures, but we were not able to find any solution of the 
corresponding equations which would produce a rotational realization. This fact supports Conjecture~\ref{NonReal}. 
The {\sc Solve} function either failed to give any 
solutions, returned complex solutions, or the solutions led to realizations where points (or lines) of different orbits 
coincide. For example, the voltage graph in Figure~\ref{fig:voltage}b gives realization shown in Figure~\ref{fig:coincide} 
where two line-orbits coincide. However, it is possible to ``perturb'' points slightly to obtain a realization in the plane which 
still resembles the $7$-cyclic symmetry, see Figure~\ref{fig:perturb}. Following Gr\"{u}nbaum \cite{Gru08}, we could try
to substitute the lines by pseudolines so that hopefully the cyclic (but not dihedral) symmetry is pre\-served.

\begin{figure}[h!]
\begin{minipage}{0.33\textwidth}
\centering
\includegraphics[scale=0.7]{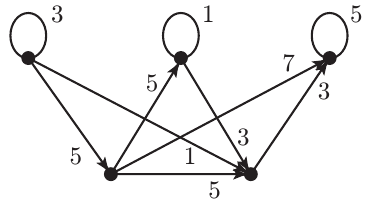}
\end{minipage}%
\begin{minipage}{0.33\textwidth}
\centering
\includegraphics[scale=0.7]{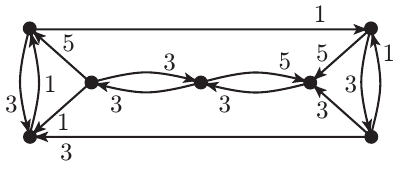}
\end{minipage}%
\begin{minipage}{0.33\textwidth}
\centering
\includegraphics[scale=0.7]{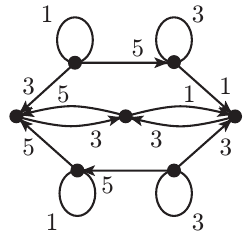}
\end{minipage}

\smallskip
\begin{minipage}{0.33\textwidth}
\centering
(a)
\end{minipage}%
\begin{minipage}{0.33\textwidth}
\centering
(b)
\end{minipage}%
\begin{minipage}{0.33\textwidth}
\centering
(c)
\end{minipage}

\caption{Voltage graphs for the Danzer graph following~\cite{census}. The Danzer graph is a $\mathbb{Z}_{14}$ covering graph over (a) 
and a $\mathbb{Z}_{10}$ covering graph over (b) and (c).}
\label{fig:PWcensus}
\end{figure}

We note in passing that we only found two bipartite voltage graphs for the Danzer graph. In \cite{census} three minimal
voltage graphs are presented: one on five vertices and two on seven vertices. They are shown in Figure~\ref{fig:PWcensus}.

We have already mentioned that the first known construction, based essentially on Theorem~\ref{PerspectiveTetrahedra}, is due to Klug~\cite{Klu}. Our aim in this Section was to find new, independent ways, especially with regard to the question of symmetry.

\section{Point-circle realizations} \label{SectionCircles}

\begin{figure}[h!]
  \begin{center}
   \includegraphics[width=0.85\textwidth]{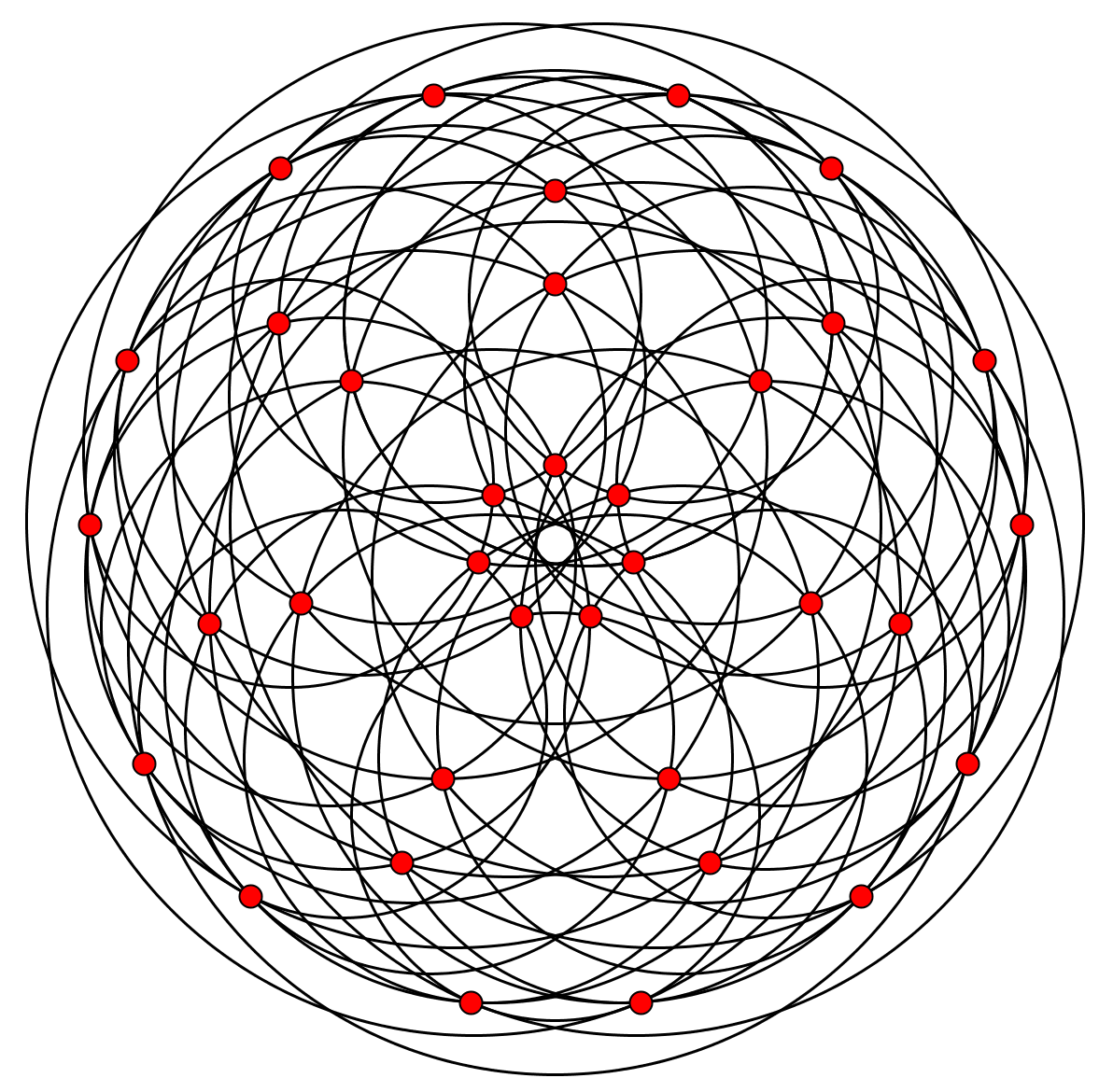}
    \caption{A point-circle realization of Danzer's configuration with 7-fold dihedral symmetry.}                
    \label{DanzerCircles}
  \end{center}
\end{figure}

The representation of the graph $O_4$ given in Figure~\ref{Odd4} has the property that each vertex-neighbourhood forms a concyclic 
set, i.e.\ a circle can be drawn through these points. Hence, the $V\!$-construction yields directly a point-circle configuration~\cite{GP}.
This configuration is shown in Figure~\ref{DanzerCircles}. Note that this realization of Danzer's configuration exhibits the dihedral 
$D_7$ symmetry; thus, this symmetry, which is the maximal possible in dimension two, can easily be achieved in this case, in contrast 
to the point-line realization (cf.\ Conjecture~\ref{NonReal}). We know that this is just a particular case of the following more general 
result~\cite{GP}.
\medskip

\begin{thm} \label{SubClifford}
For all $n\ge3$, there exists an isometric point-circle configuration of type
$$
\left(
\dbinom{2n-1}{n-1}
_{n}\,
\right).
$$
It can be obtained from the Odd graph $O_n$ by
$V\!$-construction.
\end{thm}

\noindent
In other words, for all $n\ge2$, $\DCD(n)$ can be realized as an (isometric) point-circle configuration. Moreover, it is also 
shown in~\cite{GP} that each member of this infinite series of configurations forms a subconfiguration of certain members of 
Clifford's renowned infinite series of configurations. (Here we recall that a point-circle configuration is called {\em isometric} 
if all the circles are of the same size~\cite{GP}.) Note that our example presented in Figure~\ref{DanzerCircles} is not isometric; 
however, as it is movable, its shape can be changed continuously while keeping the $D_7$ symmetry.

We expect that a point-circle realization with five-fold symmetry can also be achieved. This is supported by our following results: (1) we constructed a representation of the Danzer graph such that its Hamiltonian cycle exhibits five-fold dihedral symmetry; (2) we constructed point-circle realization of both the Steiner-Pl\"ucker and the Cayley-Salmon configuration such that both has five-fold rotational symmetry. Work in this direction is in progress.

\section* {Acknowledgements}
The authors would like to thank Primo\v{z} Poto\v{c}nik for useful discussion.

\end{document}